\documentclass[12pt]{article}
\usepackage{epsfig}
\usepackage{epsf}
\usepackage{here}
\usepackage{color}
\usepackage{pifont}

\newcommand{\DMdet}{\mbox{\boldmath ${\cal D}$}}

\newcommand{\idp}[1]{\perp \hspace{-0.25cm} \perp_{#1}}
\newcommand{\infl}[1]{\rightarrow\rightarrow_{#1}}
\newcommand{\dinf}[1]{\longrightarrow_{#1}}
\newcommand{\sli}[1]{\mbox{$\rightarrow \rightarrow\hspace{-0.60cm} / \hspace{0.40cm}_{#1}$}}
\newcommand{\wli}[1]{\mbox{$\longrightarrow \hspace{-0.60cm} / \hspace{0.40cm}_{#1}~$}}

\newcommand{\FF}{{{\cal F}}}
\newcommand{\HH}{{{\cal H}}}

\newcommand{\X}{{{\cal X}}}

\newcommand{\R}{{{\cal R}}}

\newcommand{\LL}{{\cal L}}

\newcommand{\ind}{\perp \hspace{-0.25cm} \perp}

\newcommand{\bX}{\mbox{\boldmath $X$}}

\newtheorem{Proposition}{Proposition}
\newtheorem{Definition}{Definition}

\begin{document}

\title{A general definition of influence between stochastic processes}
\author{Anne G\'egout-Petit $^{1,2}$ and Daniel Commenges$^{2,3}$  }
\maketitle
\noindent {\em 1 IMB, UMR 5251, INRIA Project CQFD, Talence, F33405, France}\\
{\em 2 Universit\'e Victor Segalen Bordeaux 2, Bordeaux, F33076, France}\\
{\em 3  INSERM, U 897, Bordeaux,  F33076, France }\\

\maketitle {\bf Summary}. \vspace{3mm}

We extend the study of weak local conditional independence (WCLI) based on a measurability condition made by Commenges and G\'egout-Petit (2009) to a larger class of processes that we call $\DMdet'$. We also give a definition related to the same concept based on certain likelihood processes, using the Girsanov theorem. Under certain conditions, the two definitions coincide on  $\DMdet'$. These results may be used in causal models in that we define what may be the largest class of processes in which influences of one component of a stochastic process on another can be described without ambiguity. From WCLI we can contruct a concept of strong local conditional independence (SCLI). When WCLI does not hold, there is a direct influence while when SCLI does not hold there is direct or indirect influence. We investigate whether WCLI and SCLI can be defined via conventional independence conditions and find that this is the case for the latter but not for the former. Finally we recall that causal interpretation does not follow from mere mathematical definitions, but requires working with a good system and with the true probability.


\section{Introduction}
The issue of causality has attracted more and more interest from statisticians in recent years.  An approach using the modelling of ``potential outcome'', often called
the counterfactual approach, has been proposed in the context of clinical trials by Rubin (1974) and further
studied by Holland (1986) among others. The counterfactual approach has been extended to the study of
longitudinal incomplete data in several papers and books (Gill and Robins, 2001; Robins et al., 2004; van der
Laan and Robins, 2002). The counterfactual approach however has been criticised (Dawid, 2000; Geneletti, 2007). Another approach directly based on dynamical models has been developed, starting with Granger (1969) and Schweder
(1970), and more recently developed using the formalism of stochastic
processes, by Aalen (1987), Florens and Foug\`ere (1996), Fosen et al. (2006) and Didelez (2007, 2008).

Recently we have given more development to the dynamical models
approach (Commenges and G\'egout-Petit, 2009) using the basic idea
of the Doob-Meyer decomposition proposed in Aalen (1987). We have
proposed a definition of weak local independence between processes
(WCLI) for a certain class of special semi-martingales (called class
$\cal{D}$) which involves the compensator of the Doob-Meyer
decomposition of the studied semi-martingale. Although it can be used in discrete time, this definition is especially adapted to continuous-time processes for which as we will see in section 4, definitions based on conventional conditional independence may fail. The aim of this paper
is to give an even more general definition of WCLI, and conversely
of direct influence. What we call direct influence of one component
$X_j$ on another component $X_k$ of a multivariate stochastic
process $\bX$ (noted $X_j \dinf {\bX} X_k $) is that $X_k$ is not
WCLI of $X_j$ (we use WCLI both as the name of the condition and as
an adjective, that is the ''I'' may mean ``independence'' or
``independent'' according to the context). This concept of influence
is a good starting point for defining \textit{causal} influence (see
section \ref{discussion}).

 In the perspective of extending WCLI to a
larger class of processes, we see two ways. The first one is to stay
in the class of semi-martingales and try to be more general about
the conditions. In particular we could use the triplet of the
characteristics of a semi-martingales. For an exact definition of
the characteristics of a semi-martingale, see Jacod and Shiryaev
(2003). Roughly speaking,  the characteristics of a semi-martingale
are represented by the triplet $(B,C , \nu)$ where $\nu$ is the
compensator of the jump part of the semi-martingale, $B$ the finite
variation part not included in $\nu$, and  $C$ is the angle bracket process of the continuous martingale. The second way is to work
with the likelihood of the process which is also tightly linked with
the characteristics of the semi-martingale. In this paper we explore these two ways, extending the WCLI definition to a
very large class of processes that we call $\DMdet '$, and showing that another definition of WCLI is possible by the use of likelihood processes. Another issue that we explore is the link between WCLI and analogous definitions based on conventional conditional independence; this angle of attack is closer to Granger (1969) proposal for time series.
The scope of the paper
is restricted to these mathematical definitions which may be useful
for discussing causality issues. In the core of the paper, we address
neither the philosophical nor the inferential issues; we discuss some of the philosophical issues in the last section.

In section 2, we recall the definition of WCLI, showing that it can
be expressed in terms of the characteristics of the
semi-martingales; this leads us to give a generalized definition of
WCLI. We also recall the definition of strong local conditional
independence (SCLI). In section 3, we propose another point of view
based on the likelihood and we show the equivalence of definitions
based on the Doob-Meyer decomposition and the property of certain
likelihood processes under certain conditions. In section 4, we show
that it is possible to define SCLI by conventional conditional independence, but
that this approach falls short for WCLI. We conclude in section 5,
where we recall the distinction between the mathematical definition
of influences and the construction of a causal interpretation.

\section{A generalization of WCLI }

\subsection{Notations and examples}

Consider a filtered space $(\Omega, \FF, (\FF_t), P)$ and a
multivariate stochastic process $\bX=(\bX_t)_{t\ge 0}$; $\bX_t$
takes  values in $\Re^m$, and the whole process $\bX$ takes values
in $D(\Re^m)$, the Skorohod space of all cadlag functions: $\Re_+
\rightarrow \Re^m$.  We suppose that all the filtrations satisfy ``the usual conditions''.
 We have $\bX=(X_j, j= 1, \ldots,m)$ where
$X_j=(X_{jt})_{t\ge 0}$. We denote by $\X_t$ the history of $\bX$ up
to time $t$, that is  $\X_t$ is the $\sigma$-field $\sigma(\bX_u,
0\le u \le t)$, and by $(\X_t)=(\X_t)_{t\ge 0}$ the families of
these histories, that is the filtration generated by $\bX$.
Similarly we shall denote by $\X_{jt}$ and $(\X_{jt})$ the histories
and the filtration associated to $X_j$. Let $\FF_t=\HH \vee \X_t$; $\HH$
may contain information known at $t=0$, in addition to the initial
value of $\bX$. We shall consider the class of special
semi-martingales in the filtration $(\FF_t$). We denote by $(B,C,\nu)$ the characteristics of the
semi-martingale $\bX$ under probability $P$, by $M_j$ the martingale
part of $X_j$, and by $M^c_j$ the continuous part of this martingale. We
denote by $(B^k,C^k,\nu^k)$ the characteristics of the semi-martingale
$X_k$ under probability $P$.

 Let us recall the definition of WCLI and see on examples how it involves the characteristics of the semi-martingale at hand. In our previous work (Commenges and G\'egout-Petit, 2009) we have imposed the two following conditions bearing on the
bracket process of the martingale $M$:

{\bf A1} $M_{j}$ and $M_{k}$ are orthogonal martingales, for all $j\ne k$;

{\bf A2} $X_j$ is either a counting process or is continuous with a deterministic bracket process, for all $j$.

We call $\DMdet$ the class of all special semi-martingales
satisfying {\bf A1} and {\bf A2}. The class of special
semi-martingales is stable by change of absolutely continuous
probability (Jacod and Shiryaev, 2003, page 43) and this is also
true for the the class $\DMdet$.

\begin{Definition}\label{WCL}[Weak conditional local independence (WCLI)]
Let $\bX$ be in the class $\DMdet$.  $X_k$ is WCLI of $X_j$ in $\bX$ on $[r,s]$ if and only if $ \Lambda
_{kt} - \Lambda_{ks}$ is $(\FF_{-jt})$-predictable on $[r,s]$, where
$\FF_{-jt}=\HH \vee \X_{-jt}$ and $\X_{-jt}=\vee _{l\ne j}\X_{-lt}$.
\end{Definition}
Generally, we assess WCLI on $[0,\tau]$, where $\tau$ is the horizon
of interest, and if $X_k$ is WCLI of $X_j$ on $[0,\tau]$, we note
$X_j \wli {\bX} X_k$; in the opposite case we say that $X_j$
directly influences $X_k$ and we note $X_j \dinf {\bX} X_k$.  A
graph representation can be given, putting a directed edge when
there is a direct influence from one node on another. If there is a
directed path from  $X_j$ to $X_k$ we say that $X_j$ influences
$X_k$ and we note $X_j \infl {\bX} X_k$.  If $X_j$ influences $X_k$
but not directly influences it, then the influence is indirect.
Inversely, if there is not directed path from $X_j$ to $X_k$,  we
say that $X_j$ does not influence $X_k$. We call this property
strong local condtional independence (SCLI), saying that  $X_j$ is
SCLI of $X_k$ and we note $X_j \sli {\bX} X_k$.

Let us see, using three examples, how the conditions {\bf A1}, {\bf
A2} and the definition of WCLI can be expressed in terms of the
characteristics of the semi-martingale $X_k$ in the filtration
$(\FF_t$).


{ \bf Example 1}: Let us consider a three-dimensional process
$\bX^3 \in \DMdet$, $\bX^3_t=(X_{1t}, X_{2t}, X_{3t})$ defined by :
\begin{equation}\label{E1}\left\{\begin{array}{ccl}
           X_{1t} & = & \int_0^t f_1(X_{1s}, X_{2s}, X_{3s})ds +  M_{1t} \\
           X_{2t} & = & \int_0^t f_2(X_{1s}, X_{2s}, X_{3s})ds +  M_{2t}  \\
           X_{3t} & = & \int_0^t f_3(X_{2s}, X_{3s})ds +  W_{3t}
         \end{array}\right.
\end{equation}
where $(M_1,M_2,W_3)$ are independent martingales and $W_3$ is a Brownian
motion. For $M_{jt}$ ($j=1,2$), we only assume that $X_{j}$ is in
the class $\DMdet$. Since $M_3=W_3$ is a
Brownian motion, $X_3$ is a diffusion, so there is no jump. In this case the characteristics of $X_3$ are
$B^3_t= \int_0^t f_3(X_{1s}, X_{2s}, X_{3s})ds$ (the finite variation process), $C^3_t= t$ (the bracket process of $W_3$) and
$\nu^3_t=0$ (the compensator of the jump part), for all $t$. From the fact that $f_3$ does not involve $X_1$, we directly see that $X_3$ is WCLI of $X_1$ in $\bX^3$. Mathematically, the ``compensator'' of the
semi-martingale $X_3$ ( called drift for a diffusion process) is equal to
$\int_0^t f_3(X_{2s}, X_{3s})ds$ for all $t$, and is thus $(\FF_{-1t})$-predictable;
 this indeed corresponds to Definition \ref{WCL} of WCLI. So, when $X_k$
is a continuous semi-martingale, the WCLI condition involves the
characteristic $B^k$ of the semi-martingale $X_k$.

{ \bf Example 2}: Let us consider the following  three-dimensional
process $\bX^3 \in \DMdet$, $\bX^3_t=(X_{1t}, X_{2t}, X_{3t})$ defined by :
\begin{equation}\label{E2}\left\{\begin{array}{ccl}
           X_{1t} & = & \int_0^t f_1(X_{1s}, X_{2s}, X_{3s})ds +  M_{1t} \\
           X_{2t} & = & \int_0^t f_2(X_{1s}, X_{2s}, X_{3s})ds +  M_{2t}  \\
           X_{3t} & = & \int_0^t \beta_3(X_{2s-}, X_{3s-})ds +  M_{3t}
         \end{array}\right.
\end{equation}
where $(M_1,M_2,M_3)$ are independent martingales and $X_3$ is a counting process.
 We do not assume the form of $M_{jt}$ ($j=1,2$). The WCLI relationships
 between the $X_i$'s are the same as in (\ref{E1}). $X_3$ is a counting process and in this
case $B^3_t=C^3_t=0$ and $\nu^3_t= \int_0^t f_3(X_{2s}, X_{3s})ds$.
 The compensator of the counting process $X_3$ is
 $\X_{-1t}$-measurable which means that $X_3$ is WCLI of $X_1$ in $\bX^3$. So, when $X_k$ is a
counting process the WCLI condition involves the characteristic $\nu^k$
of the semi-martingale $X_k$.

Thus, the WCLI condition of Commenges and G\'egout-Petit (2009)
involves  the characteristic $B$ when $X_k$ is continuous and  the
characteristic $\nu$ when $X_k$ is a counting process. In the
framework of class $\cal{D}$, condition {\bf A2} implies that these
two characteristics are never simultaneously different from zero. In
the following, we will consider processes for which $B$ and $\nu$
may be both different from zero. We consider a process each
component of which may have both a continuous and a jump part; such
a process does not belong to $\DMdet$.

{\bf Example 3}:
\begin{equation}\label{E3}\left\{\begin{array}{ccl}
X_{1t} & = & \int_0^t f_1(X_{1s},X_{2s}, X_{3s})ds + \int_0^t \sigma_{1t}dW_{1t} + \int_0^t \beta_1(X_{3s-}) dN_{1s} \\
X_{2t} & = & \int_0^t f_2(X_{1s}, X_{2s})ds + \int_0^t \sigma_{2t}dW_{2t} + \int_0^t \beta_2(X_{2s-},X_{3s-}) dN_{2s}\\
X_{3t} & = & \int_0^t f_3(X_{2s} , X_{3s})ds + \int_0^t
\sigma_{3t}dW_{3t} + \int_0^t \beta_3( X_{2s-}, X_{3s-}) dN_{3s}
         \end{array}\right.
\end{equation}
where the $W_i$'s are independent Brownian motions, the
$N_j$'s are independent Poisson Processes with intensity 1
independent of the $W_i$'s. We suppose that the $\sigma_{jt}$'s are
deterministic function of $t$, with $\sigma_{jt}>0$ $\forall t$.
It is clear that $\bX$ does not belong to class $\DMdet$.
 However, the three characteristics of the semi-martingale $X_3$ are
$B^3_t=\int_0^t f_3(X_{2s} , X_{3s})ds$, $C^3_t=\int_0^t
\sigma_{3s}ds$ and $\nu^3_t=\int_0^t \beta_3(X_{2s-} , X_{3s-})ds$.
So, $B^3_t$ and $\nu^3_t$ are$(\FF_{-1t})$-predictable: this will be the
conditions of our new WCLI available for a larger class of
semi-martingales.

\subsection{Generalized definition of WCLI}

We use the notations of the beginning of the section.
 We shall assume two conditions on  $\bX$:

{\bf A1} $M_{j}$ and $M_{k}$ are square integrable orthogonal martingales, for all
$j\ne k$.

Under assumption {\bf A1}, the jumping parts of the martingales $M_j$ and $M_k$ are orthogonal. Moreover, the characteristic $C$ of $\bX$ (the angle bracket of the continuous part of the martingale) is a
diagonal matrix. Indeed by definition of orthogonality of
semi-martingales, $C_{ij}=<M^c_i,M^c_j>=0$ for all $1 \leq i,j \leq
m$; we note $C^k=C_{kk}$.

{\bf A2'} $C^j$ is  deterministic for all $j$.

We call $\DMdet'$ the class of all special semi-martingales
satisfying {\bf A1} and {\bf A2'}. In fact, {\bf A1} and {\bf A2'} coud be merged into a single compact assuption: the  characteristic $C$ of $\bX$ is a deterministic
diagonal matrix.  $\DMdet'$ is stable by change of
absolutely continuous probability ($C$ does not change with the
probability). $\DMdet'$ is a very large class of processes: it
includes random measures, marked point processes, diffusions and diffusions with jumps.
.

\begin{Definition}[Weak conditional local independence (WCLI)]
Let $\bX$ be in the class  $\DMdet'$. $X_k$ is WCLI of $X_j$ in $\bX$ on $[r,s]$ if and only if the
characteristics $B^k$ and $\nu^k$ are such that $ B_{kt} - B_{kr}$
and $\nu_{kt} - \nu_{kr}$ are $(\FF_{-jt})$-predictable on $[r,s]$.
Equivalently we can say that $X_k$ has the same characteristic
triplet $(B^k,C^k,\nu^k)$ in $(\FF_{t})$ and in $(\FF_{-jt})$ on the
interval $[r,s]$.
\end{Definition}

This new definition coincides with that of Commenges and G\'gout-Petit (2009) for the
class $\DMdet \subset \DMdet'$.

\section{Link with the likelihood}

We consider again the three examples above with a particular
attention to the likelihood of the process $X_3$. In Example 1, we
apply Girsanov theorem to change the current probability using the
 density process $(Z_{1t}^{P/P_0})$:
\begin{equation}\label{L1}
Z^{P/P_0}_{1t} = \exp\left(\int_0^t f_3(X_{2s}, X_{3s})dX_{3s} -
\frac{1}{2}\int_0^t (f_3(X_{2s}, X_{3s}))^2ds\right).
\end{equation}

Under the assumption $E_P[\exp(\frac{1}{2}\int_0^{+\infty}
(f_3(X_{2s}, X_{3s}))^2ds)]< + \infty$, the process $Z^{P/P_0}_{1t}=
\frac{1}{Z^{P_0/P}_{1t}}$ is a $P$-martingale and the probability
$P_0$ defined by $\frac{dP_0}{dP} _{|{\FF_t}} = Z_{1t}^{P_0/P}$ for
all $t \geq 0$ is equivalent to $P$ on each $\FF_t$; moreover, under
$P_0$, $X_3$ is a Brownian motion independent of $ (M_1,M_2)$.

 In Example 2 we consider the density process $(Z^{P/P_0}_{2t})$:
\begin{equation}\label{L2}
Z^{P/P_0}_{2t} = \prod_{s \leq t}(\beta_3(X_{2s-}, X_{3s-}))^{\Delta X_{3s}}\exp\left(\int_0^t \beta_3(X_{2s-}, X_{3s-})ds \right).
\end{equation}
Under technical conditions given in L\'epingle et M\'emin (1978), it defines a new probability $P_0$ such that under $P_0$,
$X_3$ is a homogeneous Poisson process.

In Example 3, we consider the density process $(Z^{P/P_0}_{3t})$:
\begin{equation}\label{L3}
Z^{P/P_0}_{3t} = \prod_{s \leq t}(\beta_3(...))^{\Delta
X_{3s}}\exp\left(\int_0^t \frac{f_3(...)}{\sigma_{3t}}dX_{3s} +
\int_0^t (\beta_3(...)- \frac{1}{2} f^2_3(...))ds \right),
\end{equation}
where $\beta_3(...)$ stands for $\beta_3(X_{2s-}, X_{3s-})$ and $f_3(...)$ for $f_3(X_{2s}, X_{3s})$.
Under technical conditions given in L\'epingle et M\'emin (1978), $P_0$ is well defined, and
 $X_3$ is the sum of a Brownian motion with variance $\sigma^2_{3t}$ and a homogeneous Poisson process under $P_0$.

In the three cases, we see that the likelihood processes
$Z^{P/P_0}_{jt}$ are $\X_{-1t}$-measurable. That is, the
$\X_{-1t}$-measurability of the characteristics of $X_k$ implies the
$\X_{-1t}$-measurability of the likelihood process. We want to use a
measurability condition on the likelihood process for a new
definition of WCLI. We could say that  "$X_k$ is weakly locally
independent  of $X_j$ in $\bX$ if the likelihood of $X_k$ is
$\FF_{-jt}=\HH \vee \X_{-jt}$-measurable". However, we must be
cautious because the likelihood is a likelihood ratio between two
probabilities, and these probabilities give not only the
distribution of $X_k$ but that of the whole process $\bX$. So, the
reference measure $P_0$ must meet some assumptions given in the
definition of this new condition.

\begin{Definition}\label{LWCLI}[Likelihood-based weak conditional local independence (LWCLI)]\\
Let $\bX=(X_j, j= 1, \ldots,m)$ be in the class $\DMdet'$.

\begin{enumerate}
\item  Suppose the existence of a probability $P_0$  such that (i) $ P \ll P_0$, (ii)
the characteristics of the semi-martingales $X_i$'s with $i \neq k$
are the same under $P$ and $P_0$ and (iii) the $P_0$-characteristics
$(B_0^k,C_0^k,\nu_0^k)$ of the semi-martingale $X_k$ are
deterministic. We say that $X_k$ is LWCLI of $X_j$ in $\bX$ on $[0,t]$ if and only if the
likelihood ratio process $Z^{P/P_0}_{t}=\LL_{\FF_t}^{P/P_0}$ is
$(\FF_{-jt})$-measurable on $[0,t]$. We have denoted $\FF_{-jt}=\HH
\vee \X_{-jt}$ and $\X_{-jt}=\vee _{l\ne j}\X_{-lt}$.

\item $X_k$ is LWCLI of $X_j$ in $\bX$ on $[r,s]$ if and only if the  process
$\frac{Z^{P/P_0}_t}{Z^{P/P_0}_r}$ is $(\FF_{-jt})$-predictable for
all $t \in [r,s]$ for all the probabilities $P_0$ as above.
\end{enumerate}
\end{Definition}

 Let us comment the definition and the conditions
imposed to the reference probability $P_0$ in this definition in the following remarks.

\noindent{\bf Remark 1}. In the examples of this section, we have constructed $P_0$ by a change of probability. In the definition we
are in a context of likelihood writing and we suppose the existence
of a "good" reference probability.

\noindent{\bf Remark 2}. We want that the likelihood concerns $X_k$ only in a certain sense given by (ii). It was the case in the three examples considered above. (ii) is true
for instance if $<Z^{P/P_0},M^i>=0$ for all $i \neq k$. Suppose for instance that $M^k$ is not orthogonal to $M^j$ for a $j\neq k$ ( assumption {\bf A1} not true) then it is certainly not possible to find a probability which verifies (ii).

\noindent{\bf Remark 3}. We do not want
that the "relation" between $X_k$ and $X_j$ under $P$ is hidden by
the same relation under $P_0$. To make such a condition explicit, the framework of semi-martingales is again very useful. This
condition involves the characteristics $(B_0^k,C_0^k,\nu_0^k)$ of $X_k$
under $P_0$. They must be deterministic. So (iii) is linked to
assumption {\bf A2'} because $C^k$ does not change with the
probability and remains deterministic whatever the absolute
continuous change of probability. We emphasize that if {\bf A2'}
fails, that is the bracket $C^k$ is not deterministic under $P$, we
will never find a probability $P_0$ which verifies (iii). Moreover, in
the examples given above, the process $X_3$ satisfies the property
of independent increments under $P_0$. In the case of
semi-martingales, this property is verified if and only if the
triplet $(B,C,\nu)$ is deterministic under $P_0$. This is exactly
the condition (iii) of definition \ref{LWCLI}.

\noindent{\bf Remark 4}. If {\bf A1} is not satisfied, this means that at least two
components of $\bX$ have a common part of martingale: they are
driven by the same noise but we can not speak of influence of one on
the other. Condition  {\bf A2'} is different: even if $C^k$  is
driven by another component of $\bX$ we will never detect it by a
measurability condition because the characteristics $C^k$ is always
$\X_k$-measurable.

LWCLI seems to be more general than WCLI.
When $X_k$ is a diffusion with jumps (see Jacod and Shiryaev, 2003, Definition
III. 2.18), we can take for  $P_0$  the probability under which $X_k$
 is the sum of a Brownian motion and a standard
Poisson process with parameter $\lambda=1$. However, except this
standard case, the conditions required on $P_0$ are not easy to
characterize. In the good cases, we have an explicit computation of
the likelihood ratio process $Z^{P/P_0}_t$  as function of the
characteristics of $X_k$ in the probabilities $P$ and $P_0$. This
result allows us to lay down the following result:

\begin{Proposition}\label{equiv}
Suppose that  ${\bf X}$ is a $m$-dimensional diffusion with bounded
jumps process satisfying the uniqueness in law conditions  and which
belongs to the class $\DMdet'$ and suppose the existence of a
probability $P_0$ satisfying the assumptions of the definition
(\ref{LWCLI}) then $WCLI$ and $LWCLI$ are equivalent.
\end{Proposition}

{\bf Proof}: the assumptions of Proposition \ref{equiv} guarantee the explicit
computation of $Z^{P/P_0}_t$ as a function of the characteristics of
$X_k$ under $P$ and $P_0$ (Jacod and Shiryaev
2003: Theorem III. 5.19) and the uniqueness of
probability $P$ (Jacod and Shiryaev,
2003: Theorems III. 2.32 and III. 2.33) under which ${\bf X}$ has
the given characteristics.  Under these assumptions, the component
$X_k$ is of the form:
 $$dX_{kt} = f_k(t,{\bf X}_{t})dt + \sigma_k(t)
dW_{kt} + \beta_k(s, {\bf X}_{s-},z)(p(dt,dz) -
q(dt,dz)),$$
where $p(dt,dz)$ is a Poisson random measure with intensity and
$q(dt,dz)=dt \otimes F(dx)$ ($F$ is a positive $\sigma$-additive
measure on $(\R, B(\R))$. So, $B^k= \int_0^t f_k(s,{\bf X}_{s})ds$,
$\nu^k= \int_0^t\beta_k(s, {\bf X}_{s-},z)q(dt,dz)$ and $C^k=
\int_0^t \sigma^2_{3s}ds$ are the characteristics of $X_k$ under
$P$. The likelihood ratio being a function of
$(B^k,C^k,\nu^k,B_0^k,C_0^k,\nu_0^k)$, it is obvious that $WCLI$
implies $LWCLI$. Let us prove the reverse: let  $X_k$ be LWCLI of $X_j$ in $\bX$. If $B^k_t $ or
$\nu^k_t$ were not $(\FF_{-jt})$-measurable, then $Z^{P/P_0}_t$
would no longer be $(\FF_{-jt})$-measurable: this contradicts $LWCLI$ !

\section{WCLI and SCLI via conditional independence of filtrations}

Heuristically, we can state the non-influence of $X_j$ on $X_k$ by saying that, on the basis of the information at time $t$, we do not need to know $X_{ju}, u< t$ to predict  $X_k$ at $t$, or after $t$. In the previous sections, we have expressed this intuition in terms of measurability of certain processes (compensator and likelihood process). Granger (1969), working with stationary time series (in discrete time) proposed a criterion based on the variance of the prediction. Eichler and Didelez (2009) gave a clear definition of Granger non-causality in a more general setting, although still for stationary time series, and they expressed it in terms of conditional independence. They distinguish between ``strong Granger-non causality'' and ``contemporaneous independence''. With our notations, strong Granger-non causality can be expressed as:
\begin{equation} \label{SGNC} X_{ks} \idp{\FF_{-jt}} \X_{jt} , t=0, 1, \ldots; s= t+1, t+2, \ldots, t+h,\end{equation}
where $h$ is called ``horizon''.

In continuous time it is also tempting to define WCLI and SCLI in terms of conditional independence. Didelez (2008) heuristically proposed the following definition for WCLI when $\bX$ is a counting process:
\begin{equation} \label{Did} X_{kt} \idp{\FF_{-jt-}} \X_{jt-} , 0\le t\le \tau.\end{equation}
This formula attempts to express non-influence by requiring that $X_{kt}$ is independent of the past of $X_j$ given the past of the other components of $\bX$.
However as remarked in Commenges and G\'egout-Petit (2009), this condition is void in general when we consider processes in continuous time. Because conditional independence is defined via conditional probability, and in general, events of $\X_{kt}$ have conditional probabilities equal to one or zero given  $\X_{kt-}$, the condition always holds.

We now propose a rigorous definition of non-influence in continuous time based on conventional conditional independence. Moreover, since independence is defined in probability theory in terms of sigma-fields, we can state this property directly in terms of the sigma-fields $\X_{jt}, j=1,\ldots,m$, without specifying stochastic processes (as argued in Commenges, 2009, a representation of statistical models in terms of sigma-fields or filtrations is more intrinsic than in terms of random variables or stochastic processes). For simplicity we define it on $(0,\tau)$.

\begin{Definition}\label{FSCLI}Filtration-based strong conditional local independence (FSCLI)]\\
Let $(\X_{jt}), j=1,\ldots,m$ be filtrations, $\X_t=\vee_j \X_{jt}$; $\FF_t=\HH \vee \X_t$ and $\X_{-jt}=\vee _{l\ne j}\X_{-lt}$., $\FF_{-jt}=\HH \vee \X_{-jt}$.
 We say that filtration $(\X_{kt})$ is FSCLI of $(\X_{jt})$ in $\FF_t$ if and only if:
\begin{equation} \label{eFSLI} \X_{k\tau} \idp{\FF_{-jt}} \X_{jt} ,~ 0\le t\le \tau.\end{equation}
\end{Definition}

\begin{Proposition}\label{equivF}
Suppose that  ${\bf X}$ is the unique $m$-dimensional solution of a given stochastic differential equation with bounded
jumps process  and which belongs to the class $\DMdet'$, then $FSCLI$ defined on the filtrations generated by the components of $\bX$ and $SCLI$   are equivalent.
\end{Proposition}
{\bf Proof.} For a given $X_k \in \bX$, denote by $An(k)= \{l_{1k}, \ldots ,
l_{nk}\}$  the set of all the indices $l$ such that $X_l$ is
an ancestor of $X_k$. The assumptions imply that
$X_{An(k)}$ is also the unique solution a stochastic differential
equation with bounded jumps process generated by the Brownian
process ${\bf W}_{An(k)} = (W_{l_1}, \ldots ,W_{l_{nk}} )$ and the
set of orthogonal Poisson measures ${\bf P}_{An(k)}=(p_{l_1}, \ldots
,p_{l_{nk}})$.
Moreover for each $t$,  $X_{An(k)t}$ is a functional
of $({\bf W}_{An(k)s},{\bf P}_{An(k)s} , s \leq t) $. If $t \leq
\tau $, $X_{An(k)\tau}$ is a functional of $X_{An(k)t}$ and of the
processes ${\bf W}^{(t,.)}_{An(k)}$ ${\bf P}^{(t,.)}_{An(k)}$
defined by ${\bf W}^{(t,s)}_{An(k)}=({\bf W}_{An(k)s}-{\bf
W}_{An(k)t}),{\bf P}^{(s,t)}_{An(k)}=({\bf P}_{An(k)s}-{\bf
P}_{An(k)t})  , t \leq s \leq \tau) $. By the independent
 increments property of the Brownian motion and of the Poisson process, if we denote
 $\sigma_k^{(t, \tau)} = \sigma (({\bf W}^{(t,s)}_{An(k)},{\bf P}^{(t,s)}_{An(k)} , t \leq s \leq \tau)$,
 we have $\sigma_k^{(t, \tau)} \ind \FF_t$.
 Suppose that $X_j \sli {\bX} X_k$, it implies that $X_j \wli {\bX}
X_{An(k)}$ and that $X_{An(k)t}$ is $\FF_{-jt}$-measurable. Using
the previous remark and the standard properties of conditional
expectation (Jacod, Protter exercice 23.7), for $t \leq s \leq
\tau$, we have that $E[f(X_{ks})|\FF_t] = E[f(X_{ks})|\FF_{-jt}]=
G(X_{kt})$ with
$$G(x) = E[f(F(x,{\bf W}^{(t,.)}_{An(k)},{\bf P}^{(t,.)}_{An(k)}) |X_{kt}=x].$$
We have proved $SCLI \Rightarrow FSCLI$.

As for the converse, (\ref{eFSLI}) implies that $X_k$ is perfectly defined by a differential equations with jumps which does not involve
the component $X_j$ and thus $X_j \sli {\bX} X_k$.

\noindent {\bf Remark 5.} We could also define an ``horizon'' $h>0$ for FSCLI in a way analogous to formula (\ref{SGNC}).
\begin{equation} \label{FSLIhor} \X_{k,t+h} \idp{\FF_{-jt}} \X_{jt} ,~ 0\le t\le \tau-h.\end{equation}
If we make this horizon tend toward zero the FSCLI requirement tends (heuristically) to the WCLI requirement. In continuous time however, considering an infinitely small $h$ would lead to definition (\ref{Did}), which as already mentioned is void. We conclude that WCLI cannot be rigorously defined by conditional independence; we need the measurability-based definition.

\noindent {\bf Remark 6.} Didelez and Eichler (2009) also defined a concept of contemporaneous independence as:
\begin{equation} X_{ks} \idp{\FF_{t}} X_{js} , t=0, 1, \ldots; s= t+1.\end{equation}
For the same reason as for WCLI, contemporaneous independence cannot be defined in continuous time via conventional conditional independence, because $X_{kt} \idp{\FF_{t-}} X_{jt}$ in void in general. However contemporaneous independence in continuous time might be identified with the assumption of orthogonal martingales.

\noindent {\bf Remark 7.} If the time parameter is discrete, then FSLI defined by (\ref{eFSLI}) is identical to strong Granger-non causality for all horizon. Moreover FSLI defined by (\ref{FSLIhor}) for horizon $h=1$ , Granger-non causality for horizon $h=1$ and WCLI are identical.

\section{Discussion and conclusion}
{\label{discussion}
We have generalized the definition of WCLI to a larger class of processes and we have
proposed another definition through likelihood ratio processes.
Under certain conditions the two definitions are equivalent. We have also made the link with definitions based on conventional conditional independence: SCLI can be defined this way but in continuous time WCLI cannot.  These
results may be used for developing causal models. By definition, there
are direct influences where WCLI does not hold: if $X_k$ is not WCLI
of $X_j$, then $X_j$ directly influences $X_k$ in $\bX$. It is to be
noted that influence is not a simple lack of (even conditional)
independence. WCLI is a dynamical concept which differs markedly
form conventional independence concepts. Essentially because it is
dynamic, it is not symmetric, while conventional independence is. We
can have $X_k$ WCLI of $X_j$ and  $X_j$ not WCLI of $X_k$. This
provides a rich set of relationships between two components of a
stochastic process $\bX$. We have three possibilities for the
influence of $X_j$ on $X_k$:  $X_j \dinf {\bX} X_k$ (direct
influence) , $X_j \sli {\bX} X_k$ (no influence), $X_j \infl {\bX}
X_k$  and $X_j \wli {\bX} X_k$ (indirect influence). There are also
three possibilities of influence of $X_k$ on $X_j$. Thus, there are
nine possibilities for describing the relationship between two
components of a stochastic process. Of course it would be of great
interest to quantify these influences. Interesting work has been
done in this direction, in the time-series framework, by Eichler and
Didelez (2009).

The multivariate stochastic process framework is a  general
framework which incorporates a major feature of causal relationship:
time. Thus it is a natural framework to formalize causality in
statistics. It is important to know which is the most general class
of stochastic processes in which we can work for developing such a
formalisation. The class $\DMdet'$ seems to be this class.

 However this only describes a mathematical framework which is well
 suited for formalizing causality. This is why in this paper we
 speak of influence rather than causal influence. A causal
 interpretation needs an epistemological act to link the
 mathematical model to a physical reality. In particular, WCLI is
 dependent on a filtration and a probability. Commenges and
 G\'egout-Petit (2009) emphasized that the choice of the filtration
 is related to the choice of the physical system and assumed that
 there is a true probability, $P^*$, according to which the events
 of the universe are generated. Causal influences were defined as
 influences in a good (or perfect) system and under the true
 probability.

\section{References}

\setlength{\parindent}{0.0in}

\noindent Aalen O (1987) Dynamic modelling and causality. { Scand. Actuarial J.}, { 1987}: 177-190.\vspace{3mm}

\noindent Aalen O (1978) Non parametric inference for a family of
counting processes. { The Annals of Statistics} 6: 701-726.\vspace{3mm}

 Aalen O, Frigessi A (2007) What can statistics contribute
   to a causal understanding? { Scand. J. Statist.} { 34}: 155-168.\vspace{3mm}




Commenges,D. (2009) Statistical models, conventional, penalized and hierarchical likelihood. Statistics Surveys, 3: 1-17. \vspace{3mm}

\noindent Commenges D, G\'egout-Petit A (2009) A general dynamical statistical model with causal interpretation. { J. Roy. Statist. Soc. B}, 71: 1-18.\vspace{3mm}





Dawid AP (2000) Causal inference without counterfactuals. {  J. Am. Statist. Ass.} { 95}: 407-448.\vspace{3mm}



Didelez V (2007) Graphical models for composable finite Markov processes, { Scand. J. Statist.} { 34}: 169-185. \vspace{3mm}

Didelez, V (2008) Graphical models for marked point processes, { J. Roy. Statist. Soc. B} { 70}: 245-264. \vspace{3mm}

Eichler M, Didelez V (2009) On Granger causality and the effect of
interventions in time series. Research Report 09-01, Statistics Group, University of Bristol, 2009 http://www.maths.bris.ac.uk/~maxvd/\vspace{3mm}





Florens J-P, Foug\`ere D (1996) Noncausality in continuous time. { Econometrica} { 64}: 1195-1212.\vspace{3mm}

Fosen J, Ferkingstad E, Borgan O, Aalen, OO (2006) Dynamic path analysis: a new approach to analyzing time-dependent covariates. { Lifetime Data Analysis} { 12}: 143-167.  \vspace{3mm}


Geneletti S (2007) Identifying direct and indirect effects in a non-counterfactual framework. { J. Roy. Statist. Soc. A} { 69}: 199-215.\vspace{3mm}



Gill RD, Robins J (2001) Causal inference for complex longitudinal data: the continuous case. { Ann. Statist.} { 29}: 1785-1811.  \vspace{3mm}

Granger CWJ (1969) Investigating causal relations by econometric
models and cross-spectral methods. { Econometrica} { 37}: 424-438.\vspace{3mm}





Holland P. (1986) Statistics and Causal Influence. { J. Am. Statist. Ass.} { 81}: 945-960.\vspace{3mm}


Jacod J, Protter P, (2000) { Essential in probability},
Springer-Verlag.\vspace{3mm}

 Jacod J, Shiryaev
AN (2003) { Limit Theorems fo Stochastic Processes}, Berlin:
Springer-Verlag.\vspace{3mm}





L\'epingle D, M\'emin J (1978) Sur l'int\'egrabilit\'e uniforme des martingales exponentielles. Z. Wahrsch. Verw. Gebiete 42: 175--203. \vspace{3mm}







Robins J, Lok J, Gill R, van der Vaart, A (2004) Estimating the causal effect of a time-varying treatment on time-to-event using structural nested failure time models. { Statistica Neerlandica} { 58}: 271-295.\vspace{3mm}

Robin J., Scheines R, Spirtes P, Wasserman L (2003) Uniform consistency in causal inference. { Biometrika} { 90}: 491-515\vspace{3mm}

Rubin DB (1974) Estimating Causal Effects of Treatments in Randomized and Nonrandomized Studies. { Journal of Educational Psychology} { 66}: 688-
701.\vspace{3mm}




Schweder T (1970) Composable Markov processes. { J. Appl. Probab.} {7}: 400-410.\vspace{3mm}


van der Laan M, Robins J (2002) { Unified methods for censored longitudinal data and causality}, Springer, New-York. \vspace{3mm}





\end{document}